\documentclass[11pt, a4paper]{amsart}
\usepackage[top=35truemm,bottom=35truemm,left=35truemm,right=35truemm]{geometry}
\usepackage{setspace} 
\usepackage{amsmath, amsfonts,amsthm,amssymb,mathrsfs,amscd}
\usepackage{enumerate}
\usepackage[dvipdfmx]{graphicx}
\usepackage{ascmac}
\usepackage[arrow,matrix]{xy}
\usepackage{color}
\usepackage{array}
\usepackage{pifont}
\usepackage{longtable}
\usepackage{here}
\usepackage{braket}
\usepackage{mathtools}

\newtheorem{theo}{Theorem}[section]
\newtheorem*{theo*}{Theorem}
\newtheorem{defi}{Definition}[section]

\newtheorem{prop}{Proposition}[section]
\newtheorem{cor}{Corollary}[section]
\newtheorem{lem}{Lemma}[section]

\providecommand{\abs}[1]{\left\lvert#1\right\rvert}
\providecommand{\generate}[1]{\left\langle#1 \right\rangle}
\providecommand{\presentation}[2]{\left\langle#1 \middle| #2\right\rangle}

\DeclareMathOperator{\Kernel}{Ker}
\DeclareMathOperator{\Cayley}{Cay}
\DeclareMathOperator{\Gram}{Gram}

\makeatletter
\@addtoreset{equation}{section}

\makeatother

\def\disp{\displaystyle}

\def\G{\mathcal{G}}

\def\F{\mathbb{F}}

\def\N{\mathbb{N}}
\def\H{\mathbb{H}}
\def\R{\mathbb{R}}

\def\({\left(}
\def\){\right)}
\def\<{\langle}
\def\>{\rangle}

\begin{document}

\title[ ]{On the continuity of the growth rate on the space of Coxeter systems}
\author{Tomoshige Yukita}
\address{Department of Mathematics, School of Education, Waseda University, Nishi-Waseda 1-6-1, Shinjuku, Tokyo 169-8050, Japan}
\email{yshigetomo@suou.waseda.jp}
\date{}
\thanks{}

\begin{abstract} 
Floyd showed that if a sequence of compact hyperbolic Coxeter polygons converges, then so does the sequence of the growth rates of the Coxeter groups associated with the polygons \cite{Floyd}. 
For the case of the hyperbolic 3-space, Kolpakov discovered the same phenomena for specific convergent sequences of hyperbolic Coxeter polyhedra \cite{Kolpakov}. 
In this paper, we show that the growth rate is a continuous function on the space of Coxeter systems. 
This is an extension of the results due to Floyd and Kolpakov since the convergent sequences of Coxeter polyhedra discussed in \cite{Floyd, Kolpakov} give rise to that of Coxeter systems in the space of marked groups. 
\end{abstract}

\maketitle

\section{Introduction}
Let $G$ be a group and $S=(s_1,\cdots,s_n)$ be an ordered generating set of $G$. 
The pair $(G, S)$ is called a \textit{$n$-marked group}. 
For two $n$-marked groups $(G, S)$ and $(H, T)$, 
they are isomorphic if the bijection $\phi: S\to {T}$ which sends $s_i$ to $t_i$ extends to a group isomorphism from $G$ to $H$. 
The set $\G_n$ of all isomorphism classes of $n$-marked groups is called the \textit{space of $n$-marked groups}. 
$\mathcal{G}_n$ has a natural topology which makes $\mathcal{G}_n$ a compact totally disconnected space \cite{Grigorchuk}. 
The word length $\abs{x}_S$ of an element $x\in{G}$ is defined by 
\begin{equation*}
\abs{x}_S=\min{\Set{n\geq{1}|x=s_1\cdots s_n, \ s_i\in{S\cup{S^{-1}}}}}, 
\end{equation*}
and this gives rise to the word metric $d_{S}$ on $G$. 
We adopt the convention that $\abs{1_G}_S=0$ where $1_G$ is the identity element of $G$. 
We denote by $B_{(G, S)}(m)$ the ball of radius $m$ centered at $1_G$. 
Then the \textit{growth function} of $(G, S)$, denoted by $s_{(G, S)}(m)$, is defined to be the number of the elements of $B_{(G, S)}(m)$, that is, 
\begin{equation*}
s_{(G, S)}(m)=\#\Set{x\in{G}|\abs{x}_S\leq m}. 
\end{equation*}
The \textit{growth rate $\omega(G, S)$} is given by $\disp \omega(G, S)=\lim_{m\to \infty}\sqrt[m]{s_{(G, S)}(m)}$. 

\vspace{2mm}
A group $G$ is called a \textit{Coxeter group of rank $n$} if 
there exists an ordered generating set $S=(s_1,\cdots,s_n)$ such that $G$ has the following presentation: 
\begin{equation*}
G=\presentation{s_1,\cdots,s_n}{s_1^2=\cdots=s_n^2=1, \ (s_is_j)^{m_{ij}}=1 \text{ for } 1\leq i\neq j\leq n}. 
\end{equation*}
Such a generating set $S$ is called a \textit{Coxeter generating set of $G$}. 
We call the pair $(G, S)$ a \textit{Coxeter system of rank $n$}. 
The set of all Coxeter systems of rank $n$, denoted by $\mathcal{C}_n$, is called \textit{the space of Coxeter systems of rank $n$}. 
We consider the space $\mathcal{C}_n$ as a subspace of $\G_n$ (see Definition \ref{defi:3.2.}). 
In this paper, we show the following theorem (see Section \ref{sec:3}). 
\begin{theo*}
The growth rate $\omega:\mathcal{C}_n\to \R$ is a continuous function. 
\end{theo*}
In the hyperbolic geometry, Coxeter groups arise from the geometry of polyhedra as follows. 
Let $\H^d$ denote the hyperbolic $d$-space. 
A \textit{hyperbolic polyhedron} $P\subset{\H^d}$ is the intersection of finitely many closed half-spaces. 
If dihedral angles of $P$ are of the form $\pi/m \ (m\geq 2)$, then $P$ is called a \textit{hyperbolic Coxeter polyhedron}. 
The group generated by the reflections in the bounding hyperplanes of $P$, denoted by $G_P$, is a Coxeter group. 
We call $G_P$ the \textit{hyperbolic Coxeter group associated with $P$}. 
Concerning the study of the growth rates of hyperbolic Coxeter groups, we mention the results due to Floyd and Kolpakov. 
Let us write $\Delta(a_1,\cdots, a_n)$ for a hyperbolic Coxeter $n$-gon whose interior angles are $\pi/a_1,\cdots,\pi/a_n$. 
We say that $\Delta_k=\Delta(a_1(k),\cdots,a_n(k))$ converges to $\Delta=\Delta(a_1,\cdots,a_n)$ by letting $k\to \infty$ if $\disp \lim_{k\to \infty}a_i(k)=a_i$ for any $i$. 
Floyd showed that $\disp \lim_{k\to \infty}\omega(G_k)=\omega(G)$ if $\disp \lim_{k\to \infty}\Delta_k=\Delta$, where $G_k$ and $G$ are the hyperbolic Coxeter groups associated with $\Delta_k$ and $\Delta$, respectively \cite{Floyd}. 
For the case of the hyperbolic $3$-space, Kolpakov investigated a convergence of hyperbolic Coxeter polyhedra in \cite{Kolpakov}, 
and proved that if a sequence $P_k$ of hyperbolic Coxeter polyhedra converges, then so does the sequence $\omega(G_k)$ of the growth rates, 
where $G_k$ is the hyperbolic Coxeter group associated with $P_k$. 
Since the convergent sequences of hyperbolic Coxeter polygons and polyhedra give rise to that of the hyperbolic Coxeter groups in $\mathcal{C}_n$, 
our theorem can be considered as an extension of the results due to Floyd and Kolpakov (see Section \ref{sec:4}). 

\section{The space of marked groups and the growth rates}
In this section, we recall the space of marked groups (see \cite{ChampetierGuirardel, Grigorchuk} for more details). 
For readability, we give proofs for some known facts. 
\begin{defi}
For a group $G$ and its ordered generating set $S=(s_1, \cdots, s_n)$, 
the pair $(G, S)$ is called a $n$-marked group. 
Two $n$-marked groups $(G, S)$ and $(H, T)$ are said to be isomorphic if 
the bijection from $S$ to $T$ that sends $s_i$ to $t_i$ extends to a group isomorphism. 
The set of all isomorphism classes of $n$-marked groups is denoted by $\G_n$, and is called the space of $n$-marked groups. 
\end{defi}
In the sequel, we fix an ordered generating set $X=(x_1,\cdots,x_n)$ of the free group $\F_n$ of rank $n$. 
For any $n$-marked group $(G, S)$, let us denote the epimorphism from $\F_n$ onto $G$ that sends $x_i$ to $s_i$ by $\pi_{(G, S)}$. 
It is easy to see that two $n$-marked groups $(G, S)$ and $(H, T)$ are isomorphic if and only if $\Kernel{\pi_{(G, S)}}=\Kernel{\pi_{(H, T)}}$. 

\begin{defi}
Let $(G, S)$ be a $n$-marked group. 
For any $x\in{G}$, 
the word length of $x$ with respect to $S$, denoted by $\abs{x}_S$, is defined by 
\begin{equation*}
\abs{x}_S=\min{\Set{n\geq 0|x=s_1\cdots s_n, \ s_i\in{S\cup{S^{-1}}}}}. 
\end{equation*}
We denote by $B_{(G, S)}(R)$ the ball of radius $R$ centered at $1_G$, that is, 
\begin{equation*}
B_{(G, S)}(R)=\Set{x\in{G}|\abs{x}_S\leq R}. 
\end{equation*}
For abbreviation, we write $B(R)$ instead of $B_{(\F_n, X)}(R)$. 
We define a metric $d$ on $\G_n$ as follows. 
For two $n$-marked groups $(G, S)$ and $(H, T)$, 
let us denote by $v\((G, S), (H, T)\)$ the maximal radius of the balls centered at $1_{\F_n}$ where $\Kernel{\pi_{(G, S)}}=\Kernel{\pi_{(H, T)}}$, that is, 
\begin{equation*}
v\((G, S), (H, T)\)=\max{\Set{R\geq{0}|B(R)\cap{\Kernel{\pi_{(G, S)}}}=B(R)\cap{\Kernel{\pi_{(H, T)}}}}}. 
\end{equation*}
Then the metric $d$ on $\G_n$ is defined by 
\begin{equation*}
d\((G, S), (H, T)\)=e^{-v\((G, S), (H, T)\)}. 
\end{equation*}
It is known that the metric space $(\G_n, d)$ is compact. 
\end{defi}
The \textit{Cayley graph $\Cayley{(G, S)}$} of a $n$-marked group $(G, S)$ is an edge-labeled directed graph defined as follows. 
The vertex set is $G$ and two vertices $x, y$ are joined by an oriented edge from $x$ to $y$ labeled with $i$ if $x^{-1}y=s_i$. 
We consider the ball $B_{(G, S)}(R)$ as a subgraph of $\Cayley{(G, S)}$. 

\vspace{1mm}
\begin{lem}\label{lem:2.1.}
For two $n$-marked groups $(G, S), (H, T)\in{\G_n}$ and $R\geq 1$, 
the distance satisfies that $d\((G, S), (H, T)\)\leq e^{-(2R+1)}$ if and only if 
there exists an orientation and label-preserving graph isomorphism from $B_{(G, S)}(R)$ to $B_{(H, T)}(R)$ that sends the identity element of $G$ to that of $H$. 
\end{lem}
\proof 
Suppose that $d\((G, S), (H, T)\)\leq e^{-(2R+1)}$. 
By the definition of the metric on $\mathcal{G}_n$, we have that 
\begin{equation}
B(2R+1)\cap{\Kernel{\pi_{(G, S)}}}=B(2R+1)\cap{\Kernel{\pi_{(H, T)}}}. \label{eq:2.1}
\end{equation}
We define a map $\alpha:B_{(G, S)}(R)\to{B_{(H, T)}(R)}$ as follows.  
For $g=s_{i_1}^{\varepsilon_1}\cdots s_{i_k}^{\varepsilon_k} \ (k\leq R)$, 
\begin{equation*}
\alpha(s_{i_1}^{\varepsilon_1}\cdots s_{i_k}^{\varepsilon_k})=t_{i_1}^{\varepsilon_1}\cdots t_{i_k}^{\varepsilon_k}. 
\end{equation*}
First we show that the map $\alpha$ is well-defined. 
For that, suppose that an element $g\in{B_{(G, S)}(R)}$ has two expressions as follows. 
\begin{equation*}
g=s_{i_1}^{\varepsilon_1}\cdots s_{i_k}^{\varepsilon_k}=s_{j_1}^{\delta_1}\cdots s_{j_l}^{\delta_l}, 
\end{equation*}
where $0\leq k, l \leq R$. 
Then the element $w=x_{i_1}^{\varepsilon_1}\cdots x_{i_k}^{\varepsilon_k}x_{j_l}^{-\delta_l}\cdots x_{j_1}^{-\delta_1}$ of $\F_n$ is contained in $B(2R+1)\cap{\Kernel{\pi_{(G, S)}}}$. 
By the equality \eqref{eq:2.1}, we see that $\pi_{(H, T)}(w)=1_H$, and hence the map $\alpha$ is well-defined. 
It is easy to see that the map $\alpha$ is an orientation and label-preserving graph homomorphism. 
In the same manner the map $\beta:B_{(H, T)}(R)\to B_{(G, S)}(R)$ is defined by 
\begin{equation*}
\beta(t_{i_1}^{\varepsilon_1}\cdots t_{i_k}^{\varepsilon_k})=s_{i_1}^{\varepsilon_1}\cdots s_{i_k}^{\varepsilon_k}. 
\end{equation*}
Since $\beta=\alpha^{-1}$, $\alpha$ is the desired graph isomorphism. 

\begin{figure}[htbp]
\centering
\includegraphics[scale=0.5]{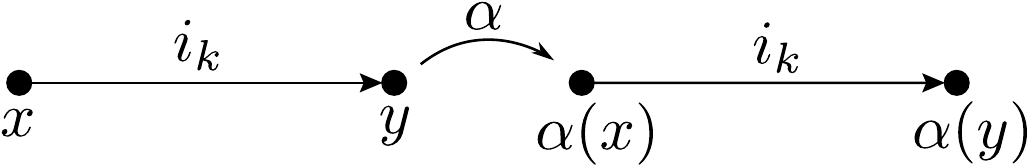}
\caption{the oriented edge from $x$ to $y$ labeled with $i_k$.}
\label{fig:figure3}
\end{figure}

\vspace{1mm}
Conversely, suppose that there exists an orientation and label-preserving graph isomorphism $\alpha:B_{(G, S)}(R)\to B_{(H, T)}(R)$ such that $\alpha(1_G)=1_H$. 
We claim that  
\begin{equation}
\alpha(s_{i_1}^{\varepsilon_1}\cdots s_{i_k}^{\varepsilon_k})=t_{i_1}^{\varepsilon_1}\cdots t_{i_k}^{\varepsilon_k}.  \label{eq:2.2}
\end{equation}
The proof is by induction on $k$. 
It is clear for the case that $k=0$. 
Set $x=s_{i_1}^{\varepsilon_1}\cdots s_{i_{k-1}}^{\varepsilon_{k-1}}$ and $y=xs_{i_k}$. 
Then the vertices $x$ and $y$ are joined by the oriented edge from $x$ to $y$ labeled with $i_k$; so are the vertices $\alpha(x)$ and $\alpha(y)$ (see  FIGURE \ref{fig:figure3}).  
By the inductive hypothesis, we have that $\alpha(x)=t_{i_1}^{\varepsilon_1}\cdots t_{i_{k-1}}^{\varepsilon_{k-1}}$. 
Since the terminal vertex of the oriented edge labeled with $i_k$ emanating from $t_{i_1}^{\varepsilon_1}\cdots t_{i_{k-1}}^{\varepsilon_{k-1}}$ is $t_{i_1}^{\varepsilon_1}\cdots t_{i_{k-1}}^{\varepsilon_{k-1}}t_{i_k}$, we have that $\alpha(y)=t_{i_1}^{\varepsilon_1}\cdots t_{i_{k-1}}^{\varepsilon_{k-1}}t_{i_k}$. 
By applying similar arguments to the case that $x=s_{i_1}^{\varepsilon_1}\cdots s_{i_{k-1}}^{\varepsilon_{k-1}}s_{i_k}^{-1}$ and $y=s_{i_1}^{\varepsilon_1}\cdots s_{i_{k-1}}^{\varepsilon_{k-1}}$, we have that $\alpha(s_{i_1}^{\varepsilon_1}\cdots s_{i_{k-1}}^{\varepsilon_{k-1}}s_{i_k}^{-1})=t_{i_1}^{\varepsilon_1}\cdots t_{i_{k-1}}^{\varepsilon_{k-1}}t_{i_k}^{-1}$. 
Therefore we obtain the equality \eqref{eq:2.2}. 
Fix an element $w=x_{i_1}^{\varepsilon_1}\cdots x_{i_k}^{\varepsilon_k}\in{B(2R+1)\cap{\Kernel{\pi_{(G, S)}}}}$. 
Let us denote by $p_w$ the closed path in $\Cayley{(G, S)}$ corresponding to $w$. 
We write $u$ for the farthest vertex of $p_w$ from $1_G$. 
If $\abs{u}_S\geq R+1$, then the length of $p_w$ must be greater than or equal to $2(R+1)$. 
Therefore the closed path $p_w$ is contained in $B_{(G, S)}(R)$. 
By the equality \eqref{eq:2.2}, the closed path of $\Cayley{(H, T)}$ corresponding to $w$ is $\alpha(p_w)$, and this implies that $B(2R+1)\cap{\Kernel{\pi_{(G, S)}}}\subset{B(2R+1)}\cap{\Kernel{\pi_{(H, T)}}}$. 
By replacing $\alpha$ with $\alpha^{-1}$, we see that $B(2R+1)\cap{\Kernel{\pi_{(H, T)}}}\subset{B(2R+1)}\cap{\Kernel{\pi_{(G, S)}}}$. \qed

\vspace{2mm}
For any $n$-marked group $(G, S)\in{\G_n}$, the \textit{growth function} $s_{(G, S)}(m)$ is defined to be the number of elements of $G$ whose lengths are at most $m$, that is, 
\begin{equation*}
s_{(G, S)}(m)=\#\Set{x\in{G}|\abs{x}_S\leq m}. 
\end{equation*}
Since $s_{(G, S)}(m)$ is submultiplicative, we have that 
\begin{equation*}
\lim_{m\to \infty}\sqrt[m]{s_{(G, S)}(m)}=\inf_{m\geq 0}\sqrt[m]{s_{(G, S)}(m)}. 
\end{equation*}
We define the \textit{growth rate}  $\omega(G, S)$ to be $\disp \lim_{m\to \infty}\sqrt[m]{s_{(G, S)}(m)}$. 
We say that a $n$-marked group $(G, S)$ has \textit{exponential growth rate} if $\omega(G, S)>1$. 

\section{The space of Coxeter systems}\label{sec:3}
In this section, we give the definitions of Coxeter groups, Coxeter matrices, and the space $\mathcal{C}_n$ of Coxeter systems of rank $n$. 
Then we shall see that $\mathcal{C}_n$ is a compact subspace of $\G_n$. 
The general reference here is \cite{BjornerBrenti}. 

\begin{defi}
Let us denote the set of the natural numbers with $\infty$ by $\widehat{\N}$. 
A $n\times{n}$ matrix $M=(m_{ij})$ over $\widehat{\N}$ is called a Coxeter matrix if 
it is symmetric and satisfies that $m_{ij}=1$ if and only if $i=j$. 
We denote the set of all Coxeter matrices of $n\times{n}$ size by $\mathcal{CM}_n$, and define a metric $D$ on $\mathcal{CM}_n$ as follows. 
\begin{equation*}
D(M, M')=\max_{1\leq i, j\leq n}\abs{\dfrac{1}{e^{m_{ij}}}-\dfrac{1}{e^{m'_{ij}}}}, 
\end{equation*}
 where $M=(m_{ij}), M'=(m'_{ij})$.
\end{defi}
Since $\mathcal{CM}_n$ is closed subset of the product space $\widehat{\N}^{n^2}$, $\mathcal{CM}_n$ is compact. 
It is easy to see that any convergent sequence $\{m(k)\}_{k\geq 1}$ of $\widehat{\N}$ is either eventually constant or converges to $\infty$, we have the following. 
\begin{lem}\label{lem:3.1.}
Let $\{M_k=\(m_{ij}(k)\)\}_{k\geq 1}$ be a sequence of Coxeter matrices. 
Then $M_k$ converges to a Coxeter matrix $M=(m_{ij})$ if and only if for any positive integer $L$, there exists $K_L\geq 0$ such that for $k\geq K_L$, 
\begin{equation*}
\begin{cases}
m_{ij}(k)\geq L & \text{ if }m_{ij}=\infty,  \\
m_{ij}(k)=m_{ij} & \text{ if }m_{ij}<\infty.
\end{cases}
\end{equation*} 
\end{lem}

\begin{defi}\label{defi:3.2.}
Any $n\times{n}$ size Coxeter matrix $M=(m_{ij})$ determines a group $G(M)$ with the presentation 
\begin{equation*}
G(M)=\presentation{s_1,\cdots,s_n}{(s_is_j)^{m_{ij}}=1\text{ for } 1\leq i\leq j\leq n}. 
\end{equation*}
A group $G$ having such a presentation is called a Coxeter group and the pair $(G, S)$ is called a Coxeter system. 
The cardinality of the generating set $S$ is called the rank of $(G, S)$. 
We denote the set of all Coxeter systems of rank $n$ by $\mathcal{C}_n\subset{\G_n}$. 
The subspace $\mathcal{C}_n$ is called the space of Coxeter systems of rank $n$. 
\end{defi}

In order to show that $\mathcal{CM}_n$ and $\mathcal{C}_n$ are homeomorphic, 
we will show that $\mathcal{C}_n$ is closed in $\G_n$. 
\begin{theo}\cite[Theorem 1.5.1, p.18]{BjornerBrenti}\label{theo:3.1.}
Let $G$ be a group and $S$ be a its generating set. 
Suppose that every element of $S$ is of order $2$. 
Then the followings are equivalent. 
\begin{itemize}
\item[(i)] The pair $(G, S)$ is a Coxeter system. 
\item[(ii)] Let $x=s_{i_1}\cdots s_{i_l}$ be a reduced word in $G$ and $s\in{S}$. 
If $\abs{sx}_S<\abs{x}_S$, then $sx=s_{i_1}\cdots \widehat{s}_{i_m}\cdots s_{i_l}$. 
\end{itemize}
\end{theo}
\begin{lem}\label{lem:3.2.}
The space of Coxeter systems $\mathcal{C}_n$ is closed in $\G_n$. 
\end{lem}
\proof Suppose that a sequence $\{(G_k, S_k)\}_{k\geq 1}$ of Coxeter systems of rank $n$ converges to a $n$-marked group $(G, S)$. 
We show that $(G, S)$ is a Coxeter system of rank $n$. 
By Lemma \ref{lem:2.1.} and the assumption that $\disp \lim_{k\to \infty}(G_k, S_k)=(G, S)$, 
for any $R\geq 0$, there exists an integer $k=k(R)\geq 1$ such that the balls $B_{(G_k, S_k)}(R)$ and $B_{(G, S)}(R)$ are isomorphic. 
By letting $R=2$, we see that every element of $S$ is of order $2$. 
Fix a reduced word $x=s_{i_1}\cdots s_{i_l}$ in $G$ and $s\in{S}$. 
Then the expression $s_{i_1}\cdots s_{i_l}$ gives rise to the shortest path from $1_G$ to $x$. 
Fix an integer $R\geq 0$ such that $x, sx\in{B_{(G, S)}(R)}$. 
Let us denote the graph isomorphism from $B_{(G, S)}(R)$ to $B_{(G_k, S_k)}(R)$ by $\alpha$. 
The assumption that $x=s_{i_1}\cdots s_{i_l}$ is a reduced word implies that $\alpha(x)=\alpha(s_{i_1})\cdots \alpha(s_{i_l})$ is a reduced word in $G_k$. 
Since $G_k$ is a Coxeter group, 
by Theorem \ref{theo:3.1.}, if $\abs{sx}_S<\abs{x}_S$, then 
\begin{equation*}
\alpha(sx)=\alpha(s)\alpha(x)=\alpha(s_{i_1})\cdots \widehat{\alpha(s_{i_m})}\cdots \alpha(s_{i_l})=\alpha(s_{i_1}\cdots \widehat{s}_{i_m}\cdots s_{i_l}). 
\end{equation*}
Therefore we see that $sx=s_{i_1}\cdots \widehat{s}_{i_m}\cdots s_{i_l}$. \qed

\begin{theo}\label{theo:3.2.}
The map $\Phi:\mathcal{CM}_n\to \mathcal{C}_n$ that sends $M$ to $G(M)$  is a homeomorphism. 
\end{theo}
\proof It is trivial that $\Phi$ is a bijection. 
If we prove that $\Phi^{-1}$ is continuous, 
the fact that $\G_n$ is compact together with Lemma \ref{lem:3.2.} implies that $\Phi^{-1}$ is a homeomorphism, since any continuous bijection from a compact space to a Hausdorff space is homeomorphism. 
Let $(G_k, S_k)$ and $(G, S)$ be Coxeter systems of rank $n$ and write $M_k=(m_{ij}(k))$ and $M=(m_{ij})$ for the Coxeter matrices corresponding to $(G_k, S_k)$ and $(G, S)$, respectively. 
Suppose that $\disp \lim_{k\to \infty}(G_k, S_k)=(G, S)$. 
Fix a defining relation $r=(s_is_j)^{m_{ij}}$ of $G$. 
The relation $r$ corresponds to the cycle of length $2m_{ij}$ in $\Cayley{(G, S)}$ labeled with $i$ and $j$ alternately. 
For the case that $m_{ij}<\infty$, 
by Lemma \ref{lem:2.1.}, 
the ball $B_{(G_k, S_k)}(m_{ij})$ must contain the cycle of length $2m_{ij}$ labeled with $i$ and $j$ alternately for sufficiently large $k$. 
Therefore $\disp \lim_{k\to \infty}m_{ij}(k)=m_{ij}$. 
Consider the case that $m_{ij}=\infty$. 
In order to obtain a contradiction, suppose that there exists an integer $R\geq 0$ such that $m_{ij}(k)\leq R$ for any $k$. 
Since $B_{(G, S)}(2R)$ does not contain the cycle labeled with $i$ and $j$ alternately, 
we see that the balls $B_{(G, S)}(2R)$ and $B_{(G_k, S_k)}(2R)$ are not isomorphic for any $k$. 
This contradicts to Lemma \ref{lem:2.1.}. 
Therefore we obtain that $\disp \lim_{k\to \infty}m_{ij}(k)=m_{ij}$ for any $i, j$, that is, the map $\Phi^{-1}$ is continuous. \qed

\begin{defi}
Let $(G, S)$ be a Coxeter system of rank $n$ and $M=(m_{ij})$ be the Coxeter matrix corresponding to $(G, S)$.  
The Gram matrix $\Gram{(G, S)}=(g_{ij})$ is a symmetric matrix of $n\times{n}$ size defined as follows. 
\begin{equation*}
g_{ij}=\begin{cases} 1 & \text{ if }i=j, \\ -\cos{\dfrac{\pi}{m_{ij}}} & \text{ if }i\neq j. \end{cases}
\end{equation*}
If $\Gram{(G, S)}$ is positive definite (resp. positive semidefinite), the Coxeter system $(G, S)$ is said to be elliptic (resp. affine) (see \cite{Vinberg} for more details). 
We call $(G, S)$ a \textit{non-affine} Coxeter system if $(G, S)$ is neither elliptic nor affine. 
\end{defi}

\begin{lem}\label{lem:3.3.}
The set of all elliptic or affine Coxeter systems of rank $n$ is closed in $\mathcal{C}_n$. 
\end{lem}
\proof 
Suppose that a sequence $\{(G_k, S_k)\}_{k\geq 1}$ of elliptic or affine Coxeter systems of rank $n$ converges to a Coxeter system $(G, S)$ of rank $n$. 
Let us write $M_k=(m_{ij}(k))$ and $M=(m_{ij})$ for the Coxeter matrices corresponding to $(G_k, S_k)$ and $(G, S)$, respectively. 
By Theorem \ref{theo:3.2.}, 
we have that $\disp \lim_{k\to \infty}g_{ij}(k)=g_{ij}$, 
where $g_{ij}(k)$ and $g_{ij}$ denote the $(i, j)$-th entries of $\Gram{(G_k, S_k)}$ and $\Gram{(G, S)}$, respectively. 
Therefore the eigenvalues $\lambda_1(k),\cdots,\lambda_n(k)$ of $\Gram{(G_k, S_k)}$ converge to the eigenvalues $\lambda_1, \cdots, \lambda_n$ of $\Gram{(G, S)}$ (see \cite[Section VI]{Bhatia}), which proves the assertion. \qed

\vspace{2mm}
The growth type of a Coxeter system $(G, S)$ is as follows. 
\begin{itemize}
\item[(i)] If $(G, S)$ is elliptic, then $G$ is finite, so that $\omega(G, S)=1$ (see \cite{Coxeter, Vinberg}). 
\item[(ii)] If $(G, S)$ is affine, then $G$ must contain a free abelian subgroup of finite index, so that $\omega(G, S)=1$ (see \cite{Coxeter, Vinberg})
\item[(iii)] If $(G, S)$ is non-affine, then $G$ must contain a free subgroup of rank at least $2$, so that $\omega(G, S)>1$ (see \cite{Harpe1987}). 
\end{itemize}

\begin{cor}\label{cor:3.1.}
A Coxeter system $(G, S)$ is non-affine if and only if $\omega(G, S)>1$. 
\end{cor}

\vspace{2mm}
There is a partial ordering $\prec$ on $\mathcal{C}_n$ defined as follows \cite{McMullen}. 
\begin{equation}
(G, S)\prec{(G', S')} \Leftrightarrow m_{ij}\leq m'_{ij} \text{ for }1\leq i,j\leq n. 
\end{equation}
Since any convergent sequence of $\widehat{\N}$ is eventually constant or converges to $\infty$, 
we obtain the following. 
\begin{lem}\label{lem:3.4.}
Let $(G_k, S_k)$ be a Coxeter system of rank $n$ that converges to $(G, S)$. 
We write $M_k=(m_{ij}(k))$ and $M=(m_{ij})$ for the Coxeter matrices corresponding to $(G_k, S_k)$ and $(G, S)$, respectively.  
Suppose that the sequence $(G_k, S_k)$ is not eventually constant. 
Then some of $m_{ij}$'s are infinity and the sequence contains an increasing subsequence. 
\end{lem}

The following theorem due to Terragni is of fundamental importance for this paper. 
For a $n$-marked group $(G, S)$, we define the function $a_{(G, S)}(m)$ to be the number of elements of $G$ whose word lengths with respect to $S$ are equal to $m$, that is, 
\begin{equation*}
a_{(G, S)}(m)=\#\Set{x\in{G}|\abs{x}_S=m}. 
\end{equation*}
We also call the function $a_{(G, S)}(m)$ the \textit{growth function} of $(G, S)$. 
\begin{theo}\cite[Theorem A., p.607]{Terragni}\label{theo:3.3.}
Let $(G, S)$ and $(G', S')$ be Coxeter systems of rank $n$. 
If $(G, S)\prec{(G', S')}$, 
then $a_{(G, S)}(m)\leq a_{(G', S')}(m)$ for $m\geq 0$. 
\end{theo}
Since $s_{(G, S)}(m)=a_{(G, S)}(0)+\cdots+a_{(G, S)}(m)$, 
if Coxeter systems $(G, S)$ and $(G', S')$ satisfy that $(G, S)\prec{(G', S')}$, 
we have that $s_{(G, S)}(m)\leq{s_{(G', S')}(m)}$ for any $m\geq 0$. 
\begin{cor}\label{cor:3.2.}
The set of all non-affine Coxeter systems of rank $n$ is closed in $\mathcal{C}_n$. 

\end{cor}
\proof 
Suppose that a sequence $\{(G_k, S_k)\}_{k\geq 1}$ of non-affine Coxeter systems of rank $n$ converges to a Coxeter system $(G, S)$. 
We show that $(G, S)$ is non-affine. 
It is trivial for the case that the sequence is eventually constant. 
For the other case, by Lemma \ref{lem:3.4.}, there exists an increasing subsequence $\{(G_{k_l}, S_{k_l})\}_{l\geq 1}$ of non-affine Coxeter systems of rank $n$. 
Theorem \ref{theo:3.3.} implies that $\omega(G, S)\geq \omega(G_{k_1}, S_{k_1})>1$, and hence $G$ is non-affine by Corollary \ref{cor:3.1.}. \qed

\vspace{1mm}
Steinberg's formula is fundamental to compute the growth rates of Coxeter systems. 
\begin{defi}
For a Coxeter system $(G, S)$ and a subset $T\subset{S}$, 
the subgroup generated by $T$, denoted by $G_T$, is called a \textit{parabolic subgroup} of $G$. 
The pair $(G_T, T)$ is a Coxeter system itself. 
We denote by $\mathcal{F}(G, S)$ the set of all the subsets of $S$ generating finite parabolic subgroups, that is, 
\begin{equation*}
\mathcal{F}(G, S)=\Set{T\subset S|\#G_T<\infty}. 
\end{equation*}
The growth series $f_{(G, S)}(z)$ is defined by 
\begin{equation*}
f_{(G, S)}(z)=\sum_{m=0}^\infty a_{(G, S)}(m)z^m. 
\end{equation*}
\end{defi}
Let us now recall the classification of finite Coxeter groups. 
In order to do that, we need the notion of \textit{Coxeter diagram}. 
\begin{defi}
For a Coxeter system $(G, S)$, 
the \textit{Coxeter diagram} $X(G, S)$ of $(G,S)$ is constructed as follows: 
its vertex set is $S$. If $m_{ij}\geq 4$, 
we join the pair of vertices by an edge and label it with $m_{ij}$. 
If $m_{ij}=3$, 
we join the pair of vertices by an edge without any label. 
For any $L\in{\widehat{\N}}$, 
we denote by $\mathcal{F}_L(G, S)$ the set of all the elements of $\mathcal{F}(G, S)$ whose Coxeter diagrams have edges labeled with $L$, that is, 
\begin{equation*}
\mathcal{F}_L(G, S)=\Set{T\in{\mathcal{F}(G, S)}| X(G_T, T) \text{ has edges labeled with }L}. 
\end{equation*}
\end{defi}
A Coxeter system $(G, S)$ is said to be \textit{irreducible} if the Coxeter diagram $X(G, S)$ is connected. 
Irreducible finite Coxeter systems of rank $n$ are classified as in TABLE \ref{tab:table1} (see \cite{BjornerBrenti}). 
For integers $m_1,\ldots,m_k\geq 1$, the polynomial $[m_1;\cdots;m_k]$ is defined by 
\begin{equation*}
[m_1;\cdots;m_k]=(1+\cdots+z^{m_1-1})\cdots (1+\cdots+z^{m_k-1}). 
\end{equation*}
The growth series of irreducible finite Coxeter systems are determined by Solomon's formula (see \cite{Solomon} for details). 
TABLE \ref{tab:table1} shows the list of the growth series of irreducible finite Coxeter systems of rank $n$. 
Note that the growth series of irreducible finite Coxeter systems are products of cyclotomic polynomials. 
\begin{table}[htbp]
\centering
\begin{tabular}{|c|c|c|}
\hline
Coxeter group&Diagram &Growth series \\
\hline
$A_n$ &
\begin{minipage}{4truecm}
\centering
\includegraphics[scale=0.5]{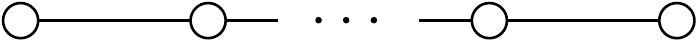}
\end{minipage}
& $[2;3;\cdots;n+1]$\\ 
\hline
$B_n$ &
\begin{minipage}{4truecm}
\centering
\includegraphics[scale=0.5]{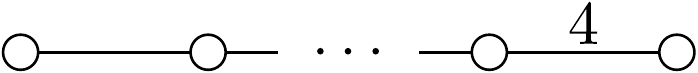}
\end{minipage}
& $[2;4;\cdots;2n]$\\ 
\hline
$D_n$ &
\begin{minipage}{4truecm}
\centering
\includegraphics[scale=0.5]{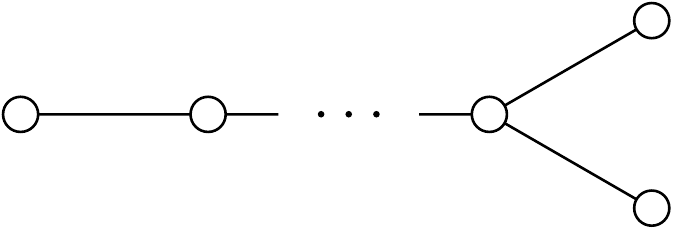}
\end{minipage}
& $[2;4;\cdots;2n-2;n]$\\ 
\hline
$E_6$ &
\begin{minipage}{4truecm}
\centering
\includegraphics[scale=0.5]{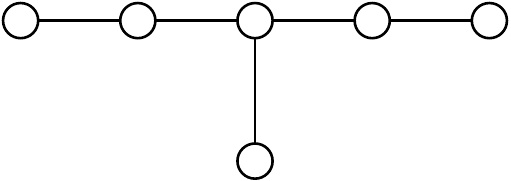}
\end{minipage}
& $[2;5;6;8;9;12]$\\ 
\hline
$E_7$ &
\begin{minipage}{4truecm}
\centering
\includegraphics[scale=0.5]{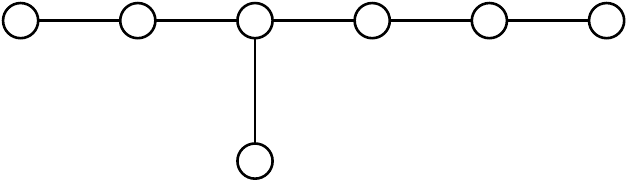}
\end{minipage}
& $[2;6;8;10;12;14;18]$\\ 
\hline
$E_8$ &
\begin{minipage}{4truecm}
\centering
\includegraphics[scale=0.5]{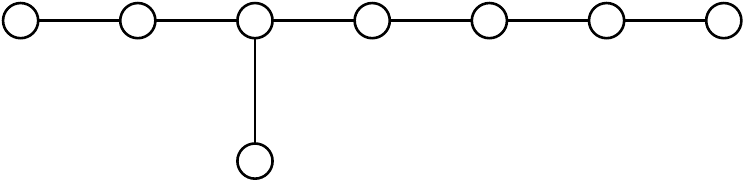}
\end{minipage}
& $[2;8;12;14;18;20;24;30]$\\ 
\hline
$F_4$ &
\begin{minipage}{4truecm}
\centering
\includegraphics[scale=0.5]{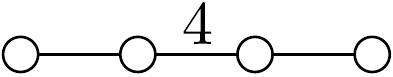}
\end{minipage}
& $[2;6;8;12]$\\ 
\hline
$H_3$ &
\begin{minipage}{4truecm}
\centering
\includegraphics[scale=0.5]{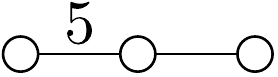}
\end{minipage}
& $[2;6;10]$\\ 
\hline
$H_4$ &
\begin{minipage}{4truecm}
\centering
\includegraphics[scale=0.5]{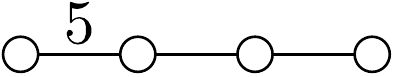}
\end{minipage}
& $[2;12;20;30]$\\ 
\hline
$I_2(m)$ &
\begin{minipage}{4truecm}
\centering
\includegraphics[scale=0.5]{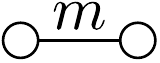}
\end{minipage}
& $[2;m]$\\ 
\hline
\end{tabular}
\caption{The classification of irreducible finite Coxeter systems of rank $n$ and their growth series. }
\label{tab:table1}
\end{table}

\vspace{1mm}
For a Coxeter system $(G, S)$, we denote the vertex sets of connected components of $X(G, S)$ by $S_1,\ldots, S_k$. 
Let $G_1,\ldots, G_k$ be the Coxeter groups generated by $S_1,\ldots,S_k$. 
The followings hold for $(G, S)$. 
\begin{equation*}
G=G_1\times{\cdots}\times{G_k}, \quad S=S_1\sqcup{\cdots}\sqcup{S_k}. 
\end{equation*} 
Then we say that $(G, S)$ is the product of $(G_1, S_1), \ldots, (G_k, S_k)$. 
\begin{lem}\label{lem:3.5.}
Let $(G, S)$ be a Coxeter system.  
If $(G, S)$ is the product of irreducible Coxeter systems $(G_1, S_1), \cdots, (G_k, S_k)$, 
then the growth series satisfies that 
\begin{equation*}
f_{(G, S)}(z)=f_{(G_1, S_1)}(z)\cdots f_{(G_k, S_k)}(z).
\end{equation*}
\end{lem}
By Lemma \ref{lem:3.5.} and the fact that the growth series of irreducible finite Coxeter systems are products of cyclotomic polynomials, we see that the growth series $f_{(G, S)}(z)$ of a finite Coxeter system $(G, S)$ has the following property. 
\begin{equation}
z^{\deg{f_{(G, S)}}}\cdot f_{(G, S)}(z^{-1})=f_{(G, S)}(z). \label{eq:3.2}
\end{equation}
The growth series of infinite Coxeter systems are written as rational functions. 
\begin{theo}[Steinberg's formula \cite{Steinberg}]
Let $G=(G, S)$ be an infinite Coxeter system. 
The following identity holds for the growth series $f_{(G, S)}(z)$. 
\begin{equation}
\dfrac{1}{f_{(G, S)}(z^{-1})}=\sum_{T\in{\mathcal{F}(G, S)}}\dfrac{(-1)^{\#T}}{f_{(G_T, T)}(z)}. \label{eq:3.3}
\end{equation}
\end{theo}
Substituting \eqref{eq:3.2} into Steinberg's formula \eqref{eq:3.3}, we obtain that 
\begin{equation}
\dfrac{1}{f_{(G, S)}(z)}=\sum_{T\in{\mathcal{F}(G, S)}} (-1)^{\#T}\dfrac{z^{d_T}}{f_{(G_T, T)}(z)}, \label{eq:3.4}
\end{equation}
where $d_T=\deg{f_{(G_T, T)}(z)}$ for $T\in{\mathcal{F}(G, S)}$. 
Set $F_{(G, S)}(z)=\dfrac{1}{f_{(G, S)}(z)}$. 
Since the growth series $f_{(G_T, T)}(z)$ is a product of cyclotomic polynomials, 
$F_{(G, S)}(z)$ is a rational function and all of those poles lie on the unit circle. 
Therefore $F_{(G, S)}(z)$ is holomorphic on the unit open disk. 
\begin{lem}\label{lem:3.6.}
Let $(G, S)$ be an infinite Coxeter system. 
The reciprocal of the growth rate $\omega(G, S)$ is the zero of $F_{(G, S)}(z)$ whose modulus is minimum among all zeros of $F_{(G, S)}(z)$. 
\end{lem}
\proof Since $G$ is infinite, we have that $\omega(G, S)=\disp \limsup_{m\to \infty}\sqrt[m]{a_{(G, S)}(m)}$ (see \cite[Section VI.C., p.182]{Harpe2000}). 
Let us denote the radius of convergence of the series $f_{(G, S)}(z)$ by $R$. 
By the Cauchy-Hadamard formula, 
\begin{equation*}
R=\dfrac{1}{\disp \limsup_{m\to \infty}\sqrt[m]{a_{(G, S)}(m)}}=\omega(G, S)^{-1}. 
\end{equation*}
Since $\disp \inf_{m\geq 0}\sqrt[m]{a_{(G, S)}(m)}=\omega(G, S)$, we obtain that $a_{(G, S)}(m)\geq{\omega(G, S)^m}$ (see \cite[56. Proposition., p.183]{Harpe2000}), and hence the series $f_{(G, S)}(\omega(G, S)^{-1})$ diverges, which proves our assertion. \qed

\begin{figure}[htbp]
\centering
\includegraphics[scale=0.5]{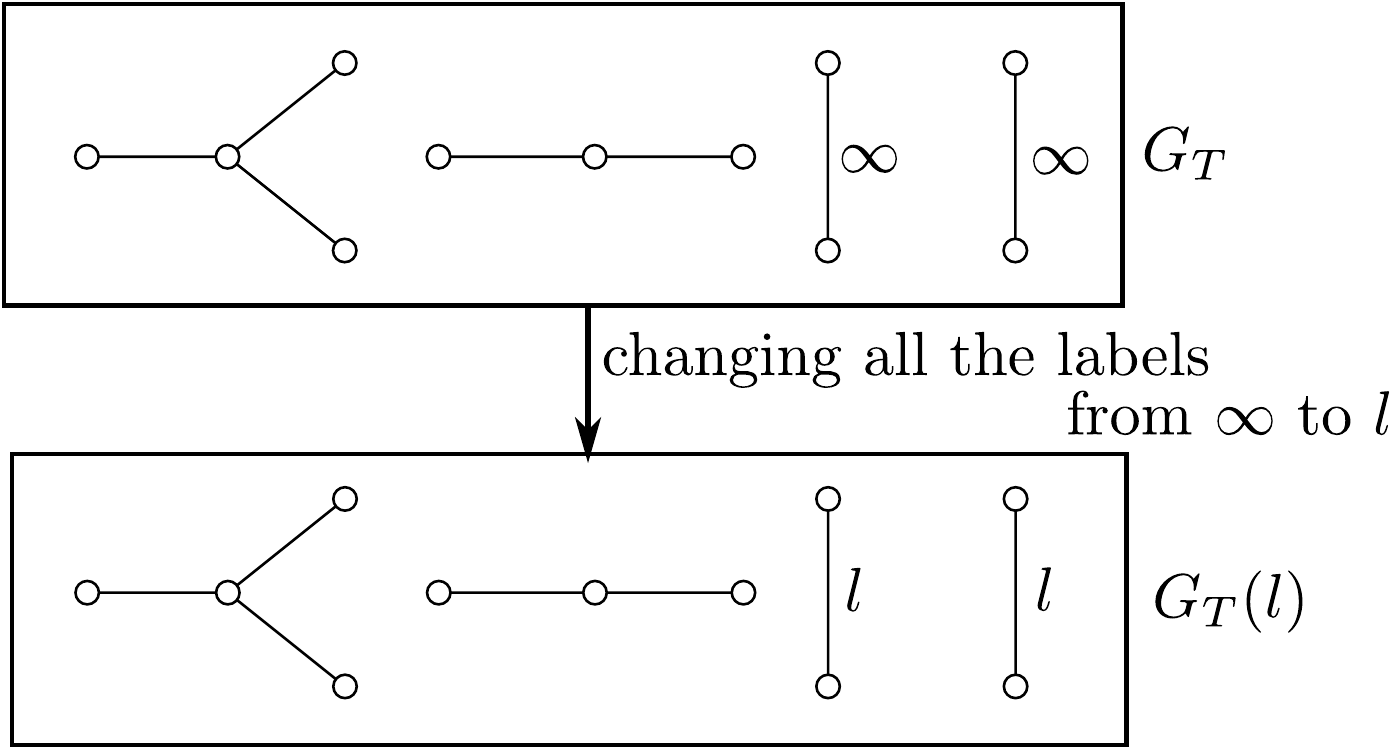}
\caption{Coxeter diagrams $X(G_T, T)$ and $X(G(l)_T, T)$.}
\label{fig:figure4}
\end{figure}

\begin{prop}\label{prop:3.1.}
Let $(G, S)$ be a non-affine Coxeter system of rank $n$ and write $M=(m_{ij})$ for the Coxeter matrix corresponding to $(G, S)$. 
For any $l\geq 6$, define the Coxeter matrix $M(l)=\(m_{ij}(l)\)$ of $n\times n$ size as follows. 
\begin{equation*}
m_{ij}(l)=\begin{cases}
l & \text{ if }m_{ij}=\infty, \\
m_{ij} & \text{ if }m_{ij}<\infty. 
\end{cases}
\end{equation*}
The Coxeter system defined by $M(l)$ is denoted by $(G(l), S(l))$. 
Then the meromorphic function $F_{(G(l), S(l))}(z)$ converges normally to $F_{(G, S)}(z)$ on the unit open disk. 
\end{prop}
\proof By Lemma \ref{lem:3.3.}, the set of non-affine Coxeter systems is open in $\mathcal{C}_n$. 
Therefore the Coxeter system $(G(l), S(l))$ is non-affine for sufficiently large $l$. 
From now on, we assume that $(G(l), S(l))$ is non-affine. 
Since the ordering on the generating set $S(l)$ of $G(l)$ is defined by the Coxeter matrix $M(l)$, 
we identify $S(l)$ with $S$ by the correspondence $s_i(l)\mapsto s_i$.  
For example, we may consider that any subset $T$ of $S$ does not only generate the parabolic subgroup $G_T$ of G but also generates the parabolic subgroup $G(l)_T$ of $G(l)$. 
For any $T\subset S$, the Coxeter diagram $X(G(l)_T, T)$ is obtained from the Coxeter diagram $X(G_T, T)$ by changing all the labels from $\infty$ to $l$ (see FIGURE \ref{fig:figure4} for an example). 
Therefore the underlying graph structure of $X(G(l)_T, T)$ is the same as that of $X(G_T, T)$. 
In this proof, we use the following notations. 
\begin{align*}
\mathcal{F}&=\mathcal{F}(G, S)=\Set{T\subset S|\#{G_T}<\infty}, \\
\mathcal{F}_L&=\mathcal{F}_L(G, S)=\Set{T\in{\mathcal{F}}|X(G_T, T)\text{ has edges labeled with }L}, \\
\mathcal{F}(l)&=\mathcal{F}(G(l), S)=\Set{T\subset S|\#{G(l)_T}<\infty}, \\
\mathcal{F}_L(l)&=\mathcal{F}_L(G(l), S)=\Set{T\in{\mathcal{F}(l)}|X(G(l)_T, T)\text{ has edges labeled with }L}, 
\end{align*}

\begin{figure}[htbp]
\centering
\includegraphics[scale=0.5]{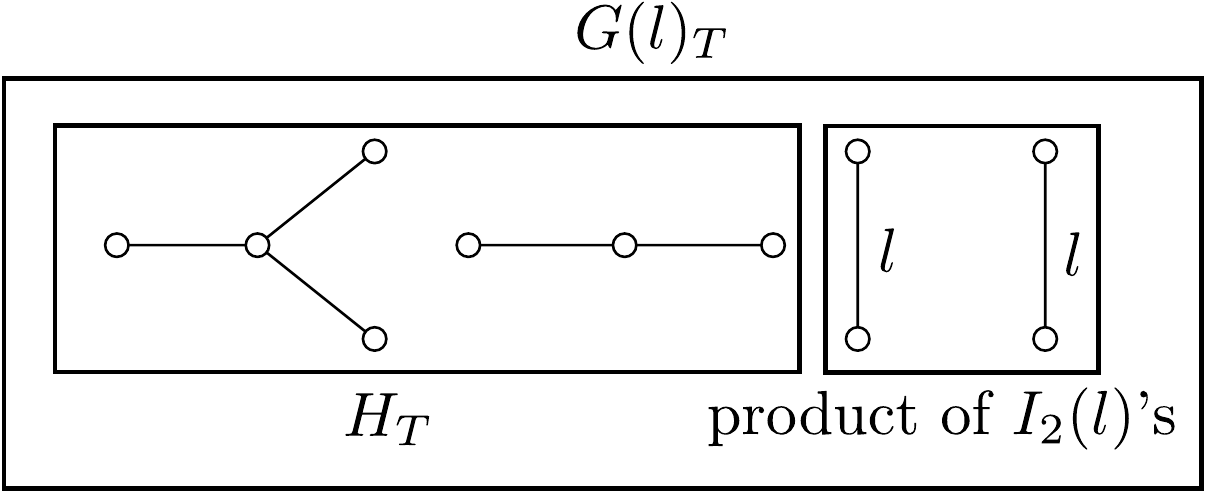}
\caption{Decomposition of $X(G(l)_T, T)$ into $H_T$ and $I_2(l)$'s.}
\label{fig:figure5}
\end{figure}

First, we show that $\mathcal{F}_l(l)=\mathcal{F}_{l'}(l')$ for any $l, l'\geq 6$. 
The classification of irreducible finite Coxeter systems implies that 
if $T\in{\mathcal{F}_l(l)}$, 
then every component of the Coxeter diagram $X(G(l)_T, T)$ having at least one edge labeled with $l$ must be $I_2(l)$. 
Hence $G(l)_T$ generated by $T\in{\mathcal{F}_l(l)}$ can be expressed as the product of the finite Coxeter group $H_T$ and $I_2(l)$'s, where the Coxeter diagram of $H_T$ has no edges labeled with $l$ (see FIGURE \ref{fig:figure5} for an example). 
Since the Coxeter diagram $X(G(l')_T, T)$ is obtained from $X(G(l)_T, T)$ by changing the labels $l$ into $l'$, $G(l')_T$ is the product of $H_T$ and $I_2(l')'s$. 
Therefore $T\in{\mathcal{F}_{l'}(l')}$ for any $T\in{\mathcal{F}_l(l)}$. 
By interchanging the roles of $l$ and $l'$, we have that $\mathcal{F}_l(l)=\mathcal{F}_{l'}(l')$. 

Given $T\in{\mathcal{F}_l(l)}$ we write $k_T$ and $I_2(l)_T$ for the number of the edges of $X(G(l)_T, T)$ labeled with $l$ and the parabolic subgroup of $G(l)_T$ generated by $T\setminus{H_T}$, respectively. 
Note that the number $k_T$ does not depend on $l\geq 6$. 
By Lemma \ref{lem:3.5.} and the fact that $(G(l)_T, T)=(H_T, T\cap{H_T})\times{(I_2(l)_T, I_2(l)_T\cap{T})}$, we have the following equality for $T\in{\mathcal{F}_l(l)}$. 
\begin{equation}
f_{(G(l)_T, T)}(z)=f_{(H_T, T\cap{H_T})}(z)f_{(I_2(l)_T, I_2(l)_T\cap{T})}(z)=f_{(H_T, T\cap{H_T})}(z)[2;l]^{k_T}. \label{eq:3.5}
\end{equation}
For $l\geq 6$ and $T\in{\mathcal{F}_l(l)}$, 
let us denote the degrees of the polynomials $f_{(G(l)_T, T)}$ and $f_{(H_T, T\cap{H_T})}$ by $d_T(l)$ and $h_T(l)$, respectively. 
By the equality \eqref{eq:3.5}, 
\begin{equation*}
d_T(l)=h_T(l)+k_Tl. 
\end{equation*}
Then the identity \eqref{eq:3.4} for $F_{(G(l), S)}(z)$ is rewritten as follows. 
\begin{align*}
F_{(G(l), S)}(z)&=\sum_{T\in{\mathcal{F}(l)}} (-1)^{\#T}\dfrac{z^{d_T(l)}}{f_{(G(l)_T, T)}(z)} \\
&=\sum_{T\in{\mathcal{F}(l)\setminus{\mathcal{F}_l(l)}}} (-1)^{\#T}\dfrac{z^{d_T(l)}}{f_{(G(l)_T, T)}(z)}+\sum_{T\in{\mathcal{F}_l(l)}} (-1)^{\#T}\dfrac{z^{d_T(l)}}{f_{(G(l)_T, T)}(z)} \\
&=\sum_{T\in{\mathcal{F}(l)\setminus{\mathcal{F}_l(l)}}} (-1)^{\#T}\dfrac{z^{d_T(l)}}{f_{(G(l)_T, T)}(z)}+\sum_{T\in{\mathcal{F}_l(l)}} (-1)^{\#T}\dfrac{z^{h_T(l)+k_Tl}}{f_{(H_T, T\cap{H_T})}(z)\cdot [2;l]^{k_T}} \\
&=\sum_{T\in{\mathcal{F}(l)\setminus{\mathcal{F}_l(l)}}}(-1)^{\#T}\dfrac{z^{d_T(l)}}{f_{(G(l)_T, T)}(z)}+\sum_{T\in{\mathcal{F}_l(l)}}(-1)^{\#T}\dfrac{z^{h_T(l)}}{f_{(H_T, T\cap{H_T})}(z)}\cdot \dfrac{z^{k_Tl}}{[2;l]^{k_T}}. 
\end{align*}
If $l$ tends to $\infty$, any parabolic subgroup $G(l)_T$ generated by $T\in{\mathcal{F}_l(l)}$ becomes an infinite Coxeter group. 
This observation implies that 
\begin{equation*}
\mathcal{F}=\mathcal{F}(l)\setminus{\mathcal{F}_{l}(l)}. 
\end{equation*}
Since the degree $d_T(l)=\deg{f_{(G(l)_T, T)}}(z)$ does not depend on $l$ for $T\in{\mathcal{F}(l)\setminus{\mathcal{F}_l(l)}}$, we obtain that 
\begin{equation*}
\sum_{T\in{\mathcal{F}(l)\setminus{\mathcal{F}_l(l)}}}(-1)^{\#T}\dfrac{z^{d_T(l)}}{f_{(G(l)_T, T)}(z)}=\sum_{T\in{\mathcal{F}}}(-1)^{\#T}\dfrac{z^{d_T}}{f_{(G(l)_T, T)}(z)}=F_{(G, S)}(z), 
\end{equation*}
and hence 
\begin{equation}
F_{(G(l), S)}(z)=F_{(G, S)}(z)+\sum_{T\in{\mathcal{F}_l(l)}}(-1)^{\#T}\dfrac{z^{h_T(l)}}{f_{(H_T, T\cap{H_T})}(z)}\cdot \dfrac{z^{k_Tl}}{[2;l]^{k_T}}. \label{eq:3.6}
\end{equation}
Now we are in a position to show that $F_{(G(l), S)}(z)$ converges normally to $F_{(G, S)}(z)$ on the unit open disk. 
Let us regard $0<\rho<1$ as fixed and write $D_\rho$ for the closed disk of radius $\rho$ centered at $0$. 
Since $f_{(H_T, T\cap{H_T})}(z)$ is a product of cyclotomic polynomials, 
there are no zeroes in $D_\rho$, and hence $\dfrac{z^{h_T(l)}}{f_{(H_T, T\cap{H_T})}(z)}$ is continuous on the compact set $D_\rho$. 
Since $\mathcal{F}_l(l)=\mathcal{F}_{l'}(l')$ for $l,l'\geq 6$, 
we may take a positive constant $M$ large enough that for any $l\geq 6, T\in{\mathcal{F}_l(l)}$ and $z\in{D_\rho}$, 
\begin{equation}
\abs{\dfrac{z^{h_T(l)}}{f_{(H_T, T\cap{H_T})}(z)}}<M. \label{eq:3.7}
\end{equation}
By multiplying $(1-z)^{k_T}$ to the denominator and numerator of $\dfrac{z^{k_Tl}}{[2;l]^{k_T}}$, we have that 
\begin{equation}
\dfrac{z^{k_Tl}}{[2;l]^{k_T}}=\(\dfrac{1-z}{1+z}\)^{k_T}\cdot \dfrac{z^{k_Tl}}{(1-z^l)^{k_T}}. \label{eq:3.8}
\end{equation}
Since the function $\dfrac{1-z}{1+z}$ is continuous on $D_\rho$, 
there exists a positive constant $M'$ such that for any $l\geq 6, T\in{\mathcal{F}_l(l)}$ and $z\in{D_\rho}$, 
\begin{equation}
\abs{\(\dfrac{1-z}{1+z}\)^{k_T}}< M'. \label{eq:3.9}
\end{equation}
By the triangle inequality,  $\abs{1-z^l}\geq 1-\rho^l$ for $l\geq 6, z\in{D_\rho}$
This observation together with the equality and inequalities \eqref{eq:3.6}, \eqref{eq:3.7}, \eqref{eq:3.8}, and \eqref{eq:3.9} give the following inequality for $z\in{D_\rho}$. 
\begin{equation*}
\abs{F_{(G(l), S)}(z)-F_{(G, S)}(z)}\leq \sum_{T\in{\mathcal{F}_l(l)}} MM' \dfrac{\rho^L}{(1-\rho^L)^{k_T}}. 
\end{equation*}
Since the cardinality of the set $\mathcal{F}_l(l)$ is constant and finite for $l\geq 6$, 
\begin{equation*}
\lim_{l\to \infty}\sup_{z\in{D_\rho}}\abs{F_{(G(l), S)}(z)-F_{(G, S)}(z)}=0, 
\end{equation*}
and this is precisely the assertion of Proposition \ref{prop:3.1.}. \qed
\begin{cor}\label{cor:3.3.}
Under the same assumption as in Proposition \ref{prop:3.1.}, 
the growth rate $\omega(G(l), S(l))$ converges to $\omega(G, S)$. 
\end{cor}
\proof Since the Coxeter system $(G(l), S(l))$ is non-affine for sufficiently large $l$, 
we may assume that $(G(l), S(l))$ is non-affine. 
In order to obtain a contradiction, suppose that $\omega(G(l), S(l))$ does not converge to $\omega(G, S)$. 
Fix $\varepsilon>0$ such that the closed disk $\overline{D(\omega(G, S)^{-1}, \varepsilon)}$ of radius $\varepsilon$ centered at $\omega(G, S)^{-1}$ does not contain $\omega(G(l), S(l))^{-1}$ for any $l$. 
Since $\omega(G, S)^{-1}$ is a zero of $F_{(G, S)}(z)$, by Proposition \ref{prop:3.1.} and Hurwitz's theorem (see \cite[Theorem, p.231]{Gamelin}), 
the disk $D(\omega(G, S)^{-1}, \varepsilon)$ contains at least one zero $z_l$ of $F_{(G(l), S(l))}(z)$ for sufficiently large $l$. 
By the triangle inequality, 
\begin{equation*}
\abs{z_l}\leq \abs{z_l-\omega(G, S)^{-1}}+\abs{\omega(G, S)^{-1}}<\varepsilon+\omega(G, S)^{-1}. 
\end{equation*}
By Theorem \ref{theo:3.3.}, $\omega(G, S)^{-1}<\omega(G(l), S(l))^{-1}$ for any $l$. 
The assumption that $\omega(G(l), S(l))^{-1}\not\in{\overline{D(\omega(G, S)^{-1}, \varepsilon)}}$ implies that 
\begin{equation*}
\omega(G(l), S(l))^{-1}>\omega(G, S)^{-1}+\varepsilon. 
\end{equation*}
By the inequalities \label{eq:3.20} and \label{eq:3.21}, 
$z_l$ is a zero of $F_{(G(l), S(l))}(z)$ whose modulus is smaller than $\omega(G(l), S(l))^{-1}$. 
This contradicts to Lemma \ref{lem:3.6.}. \qed

\begin{theo}\label{theo:3.4.}
The growth rate $\omega:\mathcal{C}_n\to \R_{\geq 1}$ is a continuous function. 
\end{theo}
\proof Let $\{(G_k, S_K)\}_{k\geq 1}$ be a convergent sequence of Coxeter systems of rank $n$ and write $(G, S)$ for the limit. 
We shall show that $\disp \lim_{k\to \infty}\omega(G_k, S_k)=\omega(G, S)$. 
The proof is divided into two cases; 
one is the case that $(G, S)$ is either elliptic or affine, and the other is the case that $(G, S)$ is non-affine. 

\vspace{2mm}
Suppose that $(G, S)$ is either elliptic or affine. 
By Corollary \ref{cor:3.2.},  the set of all elliptic or affine Coxeter systems is open in $\mathcal{C}_n$. 
This implies that $(G_k, S_k)$ is either elliptic or affine for all but finitely many $k$, and hence $\disp \lim_{k\to \infty}\omega(G_k, S_k)=1=\omega(G, S)$. 

\vspace{2mm}
Consider the case that $(G, S)$ is non-affine. 
Let us denote the Coxeter matrices corresponding to $(G_k, S_k)$ and $(G, S)$ by $M_k=(m_{ij}(k))$ and $M=(m_{ij})$, respectively. 
The assertion is trivial for the case that the sequence $(G_k, S_k)$ is eventually constant. 
We assume that $(G_k, S_k)$ is not eventually constant. 
We denote by $(G(l), S(l))$ the Coxeter system defined by the following Coxeter matrix $M(l)=(m_{ij}(l))$ for $l\geq 0$. 
\begin{equation*}
m_{ij}(l)=
\begin{cases}
l & \text{ if }m_{ij}=\infty,  \\
m_{ij} & \text{ if }m_{ij}<\infty. 
\end{cases}
\end{equation*}
By Theorem \ref{theo:3.3.} and Corollary \ref{cor:3.3.}, for any $\varepsilon>0$, there exists $L\geq 0$ such that 
\begin{equation*}
\omega(G, S)-\omega(G(L), S(L))<\varepsilon. 
\end{equation*}
By Lemma \ref{lem:3.1.} and  Theorem \ref{theo:3.2.}, 
there exists $K_L\in{\N}$ such that for $k\geq K_L$, 
\begin{equation*}
\begin{cases}
m_{ij}(k)\geq L & \text{ if }m_{ij}=\infty, \\
m_{ij}(k)=m_{ij} & \text{ if }m_{ij}<\infty. 
\end{cases}
\end{equation*}
Theorem \ref{theo:3.3.} implies that for $k\geq K_L$, 
\begin{equation*}
\omega(G, S)-\omega(G_k, S_k)\leq \omega(G, S)-\omega(G(L), S(L))<\varepsilon, 
\end{equation*}
and this proves our assertion. \qed

\section{An application to hyperbolic geometry}\label{sec:4}

\subsection{The growth rates of hyperbolic Coxeter polygons of finite volume}
We recall the result due to Floyd. 
Let $\Delta$ be a hyperbolic $n$-gon and $v_1,\cdots,v_n$ be the vertices of $\Delta$ in cyclic order. 
We call $\Delta$ a \textit{hyperbolic Coxeter $n$-gon} if the interior angles of $\Delta$ are of the form $\pi/a \ (a\geq 2)$. 
If $\Delta$ is a hyperbolic Coxeter $n$-gon, the $n$ reflections along the edges of $\Delta$ generates the Coxeter group $G(\Delta)$ of rank $n$. 
We write $\Delta(a_1,\cdots, a_n)$ for a hyperbolic Coxeter $n$-gon whose interior angle at $v_i$ equals $\pi/a_i \ (a_i\geq 2)$. 
A sequence $\{\Delta(a_1(k), \cdots, a_n(k))\}_{k\geq 1}$ of hyperbolic Coxeter $n$-gons \textit{converges} to $\Delta(a_1,\cdots, a_n)$ 
if $\disp \lim_{k\to \infty} a_i(k)=a_i$ for any $i$. 
\begin{theo}\cite[Theorem., p.476]{Floyd}\label{theo:Floyd}
Suppose that a sequence $\{\Delta(a_1(k), \cdots, a_n(k))\}_{k\geq 1}$ of hyperbolic Coxeter $n$-gons  converges to $\Delta(a_1,\cdots, a_n)$. 
Then 
\begin{equation*}
\lim_{k\to \infty}\omega(a_1(k),\cdots, a_n(k))=\omega(a_1,\cdots, a_n),
\end{equation*}
where $\omega(a_1(k),\cdots,a_n(k))$ and $\omega(a_1,\cdots,a_n)$ are the growth rates of the Coxter groups associated with $\Delta(a_1(k), \cdots, a_n(k))$ and $\Delta(a_1,\cdots,a_n)$, respectively. 
\end{theo}
We give another proof of Theorem \ref{theo:Floyd}. 
\proof Let us denote the Coxeter group associated with $\Delta(a_1,\cdots,a_n)$ by $G(a_1,\cdots,a_n)$, 
and write $M(a_1,\cdots,a_n)$ for the Coxeter matrix of $G(a_1,\cdots,a_n)$. 
By the definition of the convergence of hyperbolic Coxeter $n$-gons, 
we have that 
\begin{equation*}
\lim_{k\to \infty}M(a_1(k), \cdots, a_n(k))=M(a_1,\cdots, a_n). 
\end{equation*}
Therefore the assertion follows from Theorem \ref{theo:3.4.}. \qed

\subsection{The growth rates of hyperbolic Coxeter polyhedra of finite volume}
\begin{figure}[htbp]
\centering
\includegraphics[scale=0.5]{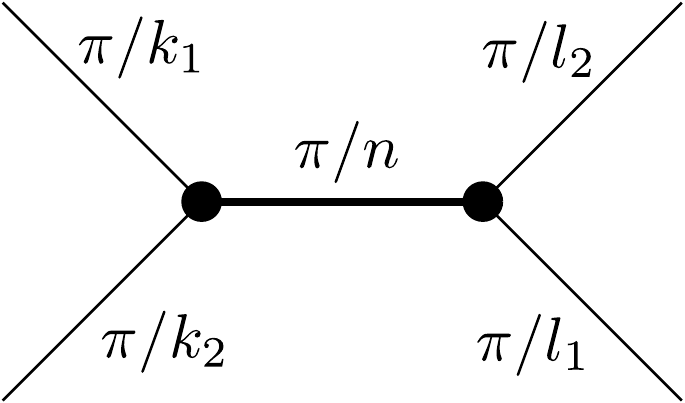}
\caption{The picture of an edge $e$ whose the endpoints are trivalent vertices}
\label{fig:figure1}
\end{figure}

\begin{figure}[htbp]
\centering
\includegraphics[scale=0.5]{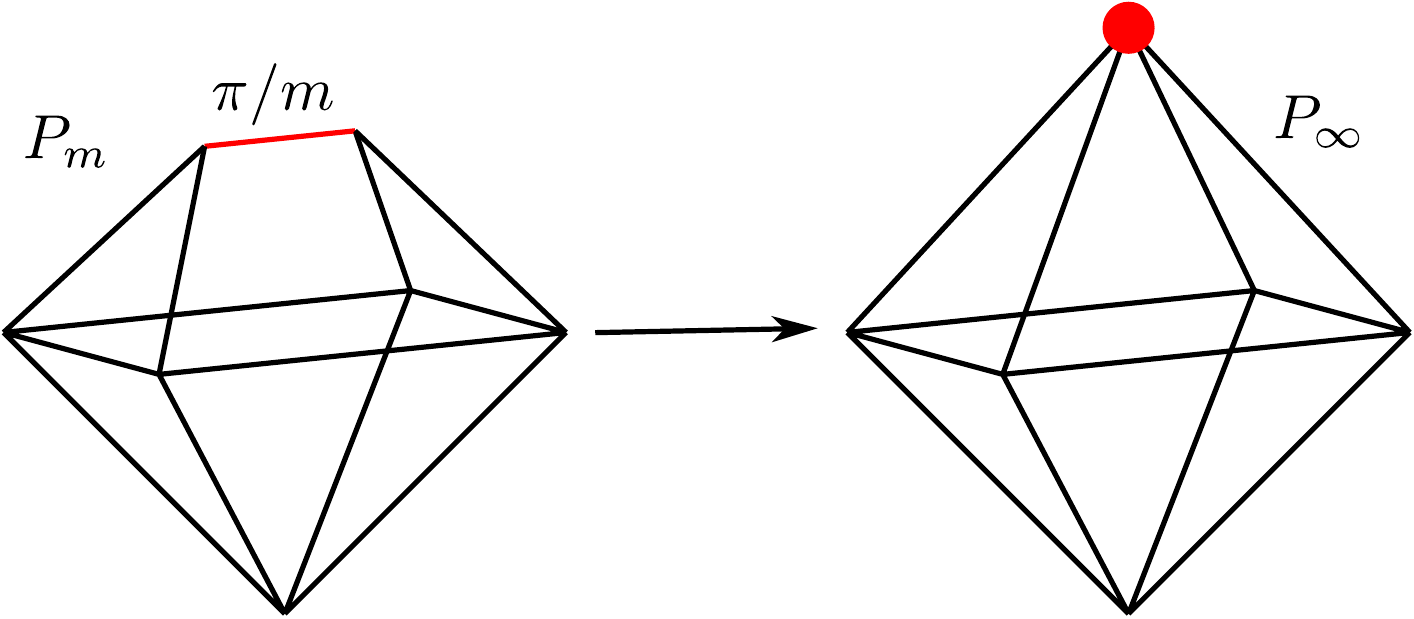}
\caption{An example of a convergent sequence of hyperbolic Coxeter polyhedra of finite-volume. The red colored edge of $P_m$ is contracted to the red colored vertex of $P_\infty$ as $m\to \infty$. }
\label{fig:figure2}
\end{figure}

We recall the edge contraction of hyperbolic Coxeter polyhedra (see \cite{Kolpakov} for more details). 
A hyperbolic Coxeter polyhedron $P$ is called a hyperbolic \textit{Coxeter} polyhedron if all dihedral angles are of the form $\pi/m \ (m\geq 2)$. 
Let $P$ be a hyperbolic Coxeter polyhedron of finite volume. 
For any edge $e$ whose the endpoints are trivalent vertices, 
we call $e$ an edge of \textit{type $\generate{k_1, k_2, n, l_1, l_2}$} if the edges incident to $e$ has dihedral angles as in FIGURE \ref{fig:figure1}. 
Suppose that a hyperbolic Coxeter polyhedron $P$ of finite volume has an edge $e$ of type $\generate{2, 2, N, 2, 2}$. 
By Andreev's theorem \cite{Andreev} and its application \cite{Kolpakov}, 
there exists a hyperbolic Coxeter polytopes $P_k$ for $N\leq k\leq \infty$ of finite volume satisfying the followings. 
\begin{itemize}
\item[(i)] For $N\leq m<\infty$, $P_m$ has the same combinatorial type as $P$. 
\item[(ii)] $P_m$ has the same dihedral angles as $P$ except for the edge $e$ and the dihedral angle $\pi/m$ at $e$. 
\item[(iii)] The combinatorial type of $P_\infty$ is obtained from $P$ by contracting the edge $e$ of $P$ to a four-valent vertex. 
\item[(iv)] $P_\infty$ has the same dihedral angles as $P$ except for the edge $e$. 
\end{itemize}
This gives us the sequence $\{P_m\}_{m\geq N}$ of hyperbolic Coxeter polyhedra converging to $P_\infty$ (see FIGURE \ref{fig:figure2}). 
We call an edge $P$ of type $\generate{2, 2, N, 2, 2}$ a \textit{contractible edge}. 

\begin{theo}\cite[Proposition 3., p.1717]{Kolpakov}\label{theo:Kolpakov}
Let $P$ be a hyperbolic Coxeter polyhedron of finite volume and $e$ be an edge of $P$ of type $\generate{2, 2, N, 2, 2}$. 
Then 
\begin{equation*}
\lim_{m\to \infty} \omega(G_m)=\omega(G_\infty), 
\end{equation*}
where $G_m$ and $G$ are the Coxeter groups associated with $P_m \ (m\geq N)$ and $P_\infty$. 
\end{theo}
We give another proof of Theorem \ref{theo:Kolpakov}. 
\proof Let us denote by $M_m$ and $M$ the Coxeter matrix of $G_m$ and $G$, repectively.  
Since the dihedral angle $\pi/m$ converges to $0$ and the other dihedral angles are constant, 
we have that 
\begin{equation*}
\lim_{m\to \infty}M_m=M. 
\end{equation*}
Therefore the assertion follows from Theorem \ref{theo:3.4.}. \qed

\subsection{The arithmetic nature of the limiting growth rates}
Let us recall Salem numbers and Pisot numbers (see \cite{Bertin} for details). 
A real algebraic integer $\alpha$ bigger than $1$ is called a \textit{Salem number} if 
$\alpha^{-1}$ is a Galois conjugate of $\alpha$ and all the other Galois conjugates lie on the unit circle. 
The set of all Salem numbers is denoted by $\mathcal{T}$. 
Parry showed that the growth rates of the Coxeter groups associated with compact hyperbolic Coxeter polygons and polyhedra are Salem numbers \cite{Parry}. 
A real algebraic integer $\alpha$ bigger than $1$ is called a \textit{Pisot number} if 
all Galois conjugates of $\alpha$ have modulus less than 1. 
The set of all Pisot numbers is denoted by $\mathcal{S}$. 
Salem numbers and Pisot numbers are closely related as follows. 
\begin{theo}\cite[Theorem 6.4.1., p.111]{Bertin}
The set $\mathcal{S}$ is contained in the closure of $\mathcal{T}$ in $\R$. 
\end{theo}
Floyd showed that the growth rates of the Coxeter groups associated with hyperbolic Coxeter polygons of finite volume are Pisot numbers \cite{Floyd}, 
and Kolpakov proved that the growth rates of the Coxeter groups associated with hyperbolic Coxeter polyhedra with single four-valent ideal vertex are Pisot numbers \cite{Kolpakov}. 
The fact that a hyperbolic Coxeter polygon and polyhedron of finite-volume are the limits of compact hyperbolic Coxeter polygons and polyhedra is of fundamental importance for their proofs. 

\begin{theo}\label{theo:4.4.}
Let $P$ be a hyperbolic Coxeter polyhedron of finite volume whose the ideal vertices are valency $4$. 
Then the growth rate $\omega(G)$ of the Coxeter group $G$ associated with $P$ is a Pisot number. 
\end{theo}
Note that Kolpakov proved that Theorem \ref{theo:4.4.} for the case that $P$ has only one ideal vertex. 
\proof  Suppose that $P$ has $N$ ideal vertices. 
By opening the ideal vertices to edges and giving sufficiently small dihedral angles $\pi/a_i \ (1\leq i\leq N)$,
we construct a compact hyperbolic Coxeter polyhedron $P(a_1,\cdots,a_N)$ whose the dihedral angles are the same as $P$ other than $\pi/a_1,\cdots,\pi/a_N$. 
By the result due to Kolpakov (see \cite[Theorem 5., p.1721]{Kolpakov}), 
the growth rate $\omega(\infty, a_2,\cdots,a_N)$ is a Pisot number. 
Since the set $\mathcal{S}$ is closed (see \cite[Theorem 6.1., p.102]{Bertin}) and $\omega(\infty, a_2, \cdots, a_N)$ converges to $\omega(\infty, \infty,a_3,\cdots,a_N)$ by Theorem \ref{theo:3.4.}, 
the growth rate $\omega(\infty, \infty, a_3,\cdots,a_N)$ is a Pisot number. 
Repeating this argument together with the equation that $\omega(G)=\omega(\infty, \cdots, \infty)$ leads to our assertion. \qed

\section*{Acknowledgement}
The author wishes to express his gratitude to Professor Koji Fujiwara for pointing out a mistake in the proof of Theorem \ref{theo:3.4.}. 
He also thanks Professor Yohei Komori for helpful discussions. 
This work was supported by JSPS Grant-in-Aid for Early-Career Scientists Grant Number JP20K14318. 

\bibliographystyle{plain}
\bibliography{reference}
\end{document}